\documentclass[a4paper,11pt]{amsart}
\usepackage{amssymb}
\usepackage{amsmath}
\usepackage{amsfonts}

\theoremstyle{plain}

\numberwithin{equation}{section}
\input{tcilatex}

\begin{document}
\title[Short Title]{On random Cameo graphs with independent edges\\
Part I: path connectivity and essential diameter }
\author{Ph. Blanchard}
\address{University of Bielefeld, Faculty of Physics, D- Bielefeld, Universit%
\"{a}tsstr.25 }
\email{Author one: blanchard@physik.uni-bielefeld.de}
\author{T. Krueger}
\email{Author~two: tkrueger@physik.uni-bielefeld.de}
\author{M. Sirugue-Collin}
\address{University of Provence, CPT, UMR 6207- 13288 Marseille Cedex 09,
France}
\email{Author three: sirugue@cpt.univ-mrs.fr}
\date{22.03.2006}

\begin{abstract}
We study growth properties of the number of paths of lenght $k$ for a
variant of Cameo graphs introduced in an earlier paper. Sharp results are
obtained for threshold for the k-path connectivity and the essential
diameter.
\end{abstract}

\maketitle

\section{Introduction}

In this paper we study phase transitions in the path-connectivity of an
inhomogeneous random graph model with power law degree distribution. The
k-path connectivity measures the average number of paths of length k between
two random chosen nodes. It has close relations to other quantities like
diameter or expected path length.

The model we use is a modification of the so-called Cameo-graphs introduced
in \cite{1}. In contrast to \cite{1} we deal here with a random graph (%
\textbf{rg}) model where all edges are independent from each other. In a
certain sense our model is an inhomogeneous extension of classical Erd\"{o}%
s-Renyi \textbf{rg} with an additional random variable $\omega $ assigned
i.i.d. to the vertices. The probability of an edge between two vertices $y$
and $x$ depends only on the $\omega -$ value of $x$ and $y$ and the vertex
set size $n$. Cameo type graphs have a nice interpretation in the context of
social network formation (see \cite{1} and \cite{2} for details).

Diameter questions for scale free graphs have been studied for different
models within the last years \cite{6,3,4,5}. The first rigorous treatment
was given by Bollobas\&Riordan \cite{3} for a precise variant (the
LCD-model) of the Albert\&Barabasi evolutionary model with preferential
attachment. In all the quoted articles a kind of branching process
approximation was used to derive the upper bounds for the diameter
respectively mean path length. In a forthcoming part II we will use the
branching process approximation for the Cameo-graphs introduced here. The
more accessible notion of path connectivity gives lower bounds on the
diameter and the expected path length and is worth studying in its own
right. The main difference between the models discussed in \cite{6,4,5} and
the Cameo type ones studied in this paper is the degree-degree correlation.
In our model the correlation is additive whereas in the mentioned models the
correlation is a multiplicative one. Both cases are interesting since they
correspond to different but somehow natural edge formation rules. The
additive case studied in this paper is in a certain sense an independent
version of the so-called $k-$ out model where each vertex generates $k$
edges independent from the other ones (with allowed multiple edges) but
still with preferences depending on the $\omega _{i}$ value of a vertex $i.$%
The degree distribution as well as evolutionary variants of this model was
studied in \cite{1}.

After completion of this article we became aware of the recent paper by
Bollobas, Janson and Riordan \cite{9} where a very general framework for
inhomogeneous random graphs with independent edges is used and in detail
analyzed. The setting therein is very similar to our approach and some of
their results can nicely be applied to the class of Cameo graphs.

\section{Definition of the model and some simple properties}

Let $\omega $ be a continuous random variable (\textbf{rv)} distributed with 
$\varphi $ such that supp$\left( \varphi \right) =\left[ 1,\infty \right) $.
Let $\Sigma $ be the set of all half infinite sequences $\hat{\omega}=\left(
\omega _{1},\omega _{2},\omega _{3},...\right) $ with i.i.d. distributed $%
\omega _{i}$ according to $\varphi $. $\Sigma $ is naturally equipped with
the product measure $\Phi $ obtained from $\varphi $. We denote by $%
C_{a,b}^{k}=\left[ \left( a_{1},b_{1}\right) ,....,\left( a_{k},b_{k}\right) %
\right] $ cylinder sets of the form $\left\{ \hat{\omega}:\omega _{i}\in %
\left[ a_{i},b_{i}\right] \text{ for }1\leq i\leq k\right\} .$ For each
element $\hat{\omega}\in \Sigma $ and fixed parameters $c>0$ and $\alpha \in 
\mathbb{R}
$ we define an associated random graph process $\left\{ \mathcal{G}\left( n,%
\hat{\omega}\right) \right\} _{n=1}^{\infty }$ as follows: $\mathcal{G}%
\left( n,\hat{\omega}\right) $ is the random graph space on vertex set $%
V_{n}=\left\{ 1;2;...;n\right\} $ with edge probability 
\begin{equation}
p_{ij}\left( n\right) :=\Pr \left\{ i\sim j\right\} =\min \left\{ \frac{c}{%
\sum_{k=1}^{n}\frac{1}{\varphi ^{\alpha }\left( \omega _{k}\right) }}\left( 
\frac{1}{\varphi ^{\alpha }\left( \omega _{i}\right) }+\frac{1}{\varphi
^{\alpha }\left( \omega _{j}\right) }\right) ,1\right\}  \label{1}
\end{equation}%
where all edges are drawn independent of each other. As will be shown later
the parameter $c$ determines the expected edge density like in classical Erd%
\"{o}s\&Renyi random graphs and the affinity parameter $\alpha $ directly
relates to the shape of the degree distribution. The value of the \textbf{rv}
$\omega _{i}$ is called the weight of vertex $i$. An element $G\in \mathcal{G%
}\left( n,\hat{\omega}\right) $ with edge set $E:=E\left( G\right) $ has
therefore probability $\Pr \left( G\right) =\prod\limits_{\left( i,j\right)
\in E}p_{ij}\prod\limits_{\left( k,l\right) \notin E}\left( 1-p_{kl}\right) $
(if no confusion arises we will drop the $n$ - dependence in $p_{ij}\left(
n\right) $).

We are interested in the asymptotic properties of $\mathcal{G}\left( n,\hat{%
\omega}\right) $ as $n\rightarrow \infty $ for $\Phi -$ typical realizations 
$\hat{\omega}\in \Sigma $. $\hat{\omega}$ is called typical if for any
cylinder set $C_{a,b}^{k}$ with $\Phi \left( C_{a,b}^{k}\right) >0$ one has 
\begin{equation}
\lim\limits_{n\rightarrow \infty }\frac{1}{n}\sum\limits_{i=0}^{n-1}\mathbf{1%
}_{C_{a,b}^{k}}\left( \sigma ^{i}\hat{\omega}\right) =\Phi \left(
C_{a,b}^{k}\right)
\end{equation}%
where $\sigma $ is the usual left-shift and $\mathbf{1}_{A}\left( \hat{\omega%
}\right) $ is the indicator function for the event $\hat{\omega}\in A\subset
\Sigma $. We say that a property $\emph{P}$ holds with high probability (%
\textbf{whp}) for $\left\{ \mathcal{G}\left( n,\hat{\omega}\right) \right\} $
if $\lim\limits_{n\rightarrow \infty }\Pr \left( \emph{P}\text{ holds for }%
G\in \mathcal{G}\left( n,\hat{\omega}\right) \right) =1$. Note that $\left\{ 
\mathcal{G}\left( n,\hat{\omega}\right) \right\} _{n=1}^{\infty }$ admits a
kind of natural filtration such that $\mathcal{G}\left( n+1,\hat{\omega}%
\right) $ is essentially obtained from $\mathcal{G}\left( n,\hat{\omega}%
\right) $ by adding vertex $n+1$ with weight $\omega _{n+1}$ and creating
new edges to vertex $n+1$ according to formula \ref{1} (where the sum in the
normalization runs now from $1$ to $n+1$) and eliminate existing edges in $%
\mathcal{G}\left( n,\hat{\omega}\right) $ with probability $p_{ij}\left(
n\right) -p_{ij}\left( n+1\right) $.

The expected degree $\bar{d}\left( i\right) $ of vertex $i\leq n$ in $%
\mathcal{G}\left( n,\hat{\omega}\right) $ is given by 
\begin{eqnarray}
\bar{d}\left( i\right) &=&\sum\limits_{j\in V_{n};j\neq i}\frac{c}{%
\sum_{k=1}^{n}\frac{1}{\varphi ^{\alpha }\left( \omega _{k}\right) }}\left( 
\frac{1}{\varphi ^{\alpha }\left( \omega _{i}\right) }+\frac{1}{\varphi
^{\alpha }\left( \omega _{j}\right) }\right)  \notag \\
&=&c+\frac{\left( n-2\right) }{\sum_{k=1}^{n}\frac{1}{\varphi ^{\alpha
}\left( \omega _{k}\right) }}\frac{c}{\varphi ^{\alpha }\left( \omega
_{i}\right) }  \notag \\
&=&c\left( 1+\frac{\left( n-2\right) }{\varphi ^{\alpha }\left( \omega
_{i}\right) \sum_{k=1}^{n}\frac{1}{\varphi ^{\alpha }\left( \omega
_{k}\right) }}\right)  \label{b}
\end{eqnarray}%
and the expected number of edges $E_{n}\left( \hat{\omega}\right) $ by 
\begin{eqnarray}
E_{n}\left( \hat{\omega}\right) &=&\frac{1}{2}\sum\limits_{i,j\in
V_{n};j\neq i}\frac{c}{\sum_{k=1}^{n}\frac{1}{\varphi ^{\alpha }\left(
\omega _{k}\right) }}\left( \frac{1}{\varphi ^{\alpha }\left( \omega
_{i}\right) }+\frac{1}{\varphi ^{\alpha }\left( \omega _{j}\right) }\right) 
\notag \\
&=&\frac{1}{2}\sum_{i}c\left( 1+\frac{\left( n-2\right) }{\varphi ^{\alpha
}\left( \omega _{i}\right) \sum_{k=1}^{n}\frac{1}{\varphi ^{\alpha }\left(
\omega _{k}\right) }}\right)  \notag \\
&=&\frac{nc}{2}+\frac{c\left( n-2\right) }{2\sum_{k=1}^{n}\frac{1}{\varphi
^{\alpha }\left( \omega _{k}\right) }}\sum_{i}\frac{1}{\varphi ^{\alpha
}\left( \omega _{i}\right) }=c\left( n-1\right)  \label{a}
\end{eqnarray}

By the ergodic theorem we have $\Phi -$ almost surely 
\begin{equation}
A\left( \alpha \right) :=\lim\limits_{n\rightarrow \infty }\frac{1}{n}%
\sum_{i=1}^{n}\frac{1}{\varphi ^{\alpha }\left( \omega _{i}\right) }%
=\int_{1}^{\infty }\varphi ^{1-\alpha }\left( \omega \right) d\omega
\end{equation}

Let $B\left( \varphi \right) $ be the set of $\alpha -$ values for which $%
A\left( \alpha \right) <\infty $ (in case $\varphi $ decays faster then any
polynomial, $B\left( \varphi \right) $ is just the set $\left( -\infty
,1\right) $). For $\alpha \in B\left( \varphi \right) $ we get the following
estimate on the asymptotic value of $\bar{d}\left( i\right) :=\mathbb{E}_{%
\mathcal{G}\left( n,\hat{\omega}\right) }\left[ d\left( i\right) \right] $
for typical $\hat{\omega}$:%
\begin{equation}
\lim\limits_{n\rightarrow \infty }\bar{d}\left( i\right) =c\left( 1+\frac{1}{%
\varphi ^{\alpha }\left( \omega _{i}\right) A\left( \alpha \right) }\right)
\end{equation}%
Note that the expected value of the number of edges in $\mathcal{G}\left( n,%
\hat{\omega}\right) $ does not depend on $\hat{\omega}$. Since the last
formula gives the expected degree of a vertex $x$ conditioned to $\omega
_{i} $ we can estimate the expected largest degree value in $\mathcal{G}%
\left( n,\hat{\omega}\right) $ for typical $\hat{\omega}$ by estimating the
asymptotic of $f\left( n,\hat{\omega}\right) :=\max\limits_{i\leq n}\hat{%
\omega}$.

Defining $F\left( z\right) :=\int\limits_{1}^{z}\varphi \left( \omega
\right) d\omega $ and $F^{\ast }\left( z\right) =1-F\left( z\right) $ we
have the following estimation on $f\left( n,\hat{\omega}\right) :$

\textbf{Lemma 1: }$f\left( n,\hat{\omega}\right) =\left[ F^{\ast }\right]
^{-1}\left( \frac{1}{n^{1+o\left( 1\right) }}\right) =F^{-1}\left( 1-\frac{1%
}{n^{1+o\left( 1\right) }}\right) ,$ $\Phi -$ a.s.

\textbf{Proof:} Clearly for any $\varepsilon >0$ one has%
\begin{eqnarray}
\Pr \left\{ f\left( n,\hat{\omega}\right) >\left[ F^{\ast }\right]
^{-1}\left( \frac{1}{n^{1+\varepsilon }}\right) \right\}
&=&1-\prod\limits_{i=1}^{n}\Pr \left\{ \omega _{i}\leq \left[ F^{\ast }%
\right] ^{-1}\left( \frac{1}{n^{1+\varepsilon }}\right) \right\}  \notag \\
&=&1-\left( 1-\frac{1}{n^{1+\varepsilon }}\right) ^{n}\underset{n\rightarrow
\infty }{\rightarrow }0
\end{eqnarray}%
and 
\begin{eqnarray}
\Pr \left\{ f\left( n,\hat{\omega}\right) <\left[ F^{\ast }\right]
^{-1}\left( \frac{1}{n^{1-\varepsilon }}\right) \right\}
&=&\prod\limits_{i=1}^{n}\Pr \left\{ \omega _{i}\leq \left[ F^{\ast }\right]
^{-1}\left( \frac{1}{n^{1-\varepsilon }}\right) \right\}  \notag \\
&=&\left( 1-\frac{1}{n^{1-\varepsilon }}\right) ^{n}\underset{n\rightarrow
\infty }{\rightarrow }0
\end{eqnarray}

Since the sequences $\left\{ 1-\left( 1-\frac{1}{n^{1+\varepsilon }}\right)
^{n}\right\} $ and $\left\{ \left( 1-\frac{1}{n^{1-\varepsilon }}\right)
^{n}\right\} $ are both summable for any $\varepsilon >0$ the Borel-Cantelli
lemma implies the pointwise statement in lemma 1.$\square $

In the next paragraph we will need some estimations of the sums $\frac{1}{n}%
\sum_{i=1}^{n}\frac{1}{\varphi ^{\beta }\left( \omega _{i}\right) }$ for $%
1<\beta \notin B\left( \varphi \right) $. Here one cannot apply directly the
ergodic theorem since $\frac{1}{\varphi ^{\beta }\left( \omega \right) }%
\notin L_{\Phi }^{1}.$ But under mild assumptions on the monotonicity of $%
\varphi $ one can show that the derived \textbf{rv} $\frac{1}{\varphi
^{\beta }\left( \omega \right) }$ is power law like distributed (lemma 2)
which in turn can be used to estimate the above mentioned ergodic sum.

\textbf{Lemma 2: }

\textbf{i) }Let $\varphi \in C^{2}$and $D^{2}\left( \varphi ^{\mu }\right)
\neq 0$ for $\left\vert \mu \right\vert \in \left( 0,1\right] $ and $\omega
>\omega _{0}\left( \mu \right) $ (this implies also that $\varphi $ is
faster decaying then any power law). Then the distribution $\psi \left(
y\right) $ of $y:=\left[ \varphi \left( \omega \right) \right] ^{-\beta
}:=\varphi \left( \omega \right) ^{-\beta }$ has asymptotic density $\psi
\left( y\right) =\frac{1}{y^{1+\frac{1}{\beta }+o_{y}\left( 1\right) }}$.

\textbf{ii)} In case $\varphi \left( \omega \right) =\frac{const}{\omega
^{\gamma }}$ the distribution of $y:=\left[ \varphi \left( \omega \right) %
\right] ^{-\beta }$ is given by $\psi \left( y\right) =\frac{1}{y^{1+\frac{1%
}{\beta }-\frac{1}{\beta \gamma }+o_{y}\left( 1\right) }}$.

For a proof see theorem 2 in \cite{2}. For convenience of the reader we give
a sketch of the proof in the appendix.

Together with lemma 1 we obtain in lemma 3 the following estimate :

\textbf{Lemma 3}: Under the assumptions of lemma 2 and for $1<\beta \notin
B\left( \varphi \right) $ and $\Phi $ a.s. the following properties hold:

\textbf{i)}%
\begin{equation}
\frac{\int_{1}^{F^{-1}\left( 1-\frac{1}{n}\right) }\varphi ^{1-\beta }\left(
\omega \right) d\omega }{\frac{1}{n}\sum_{i=1}^{n}\frac{1}{\varphi ^{\beta
}\left( \omega _{i}\right) }}=n^{o\left( 1\right) }  \label{aa}
\end{equation}

\textbf{ii)} under the assumptions of lemma 2i one has 
\begin{equation}
\int_{1}^{F^{-1}\left( 1-\frac{1}{n}\right) }\varphi ^{1-\beta }\left(
\omega \right) d\omega =\int_{y_{\min }}^{n^{\beta +o_{n}\left( 1\right)
}}y^{-\frac{1}{\beta }+o_{y}\left( 1\right) }dy=n^{\beta -1+o_{n}\left(
1\right) }  \label{ab}
\end{equation}

\textbf{iii)} for $\varphi \left( \omega \right) =\frac{const}{\omega
^{\gamma }}$ one has 
\begin{equation}
\int_{1}^{F^{-1}\left( 1-\frac{1}{n}\right) }\varphi ^{1-\beta }\left(
\omega \right) d\omega =\int_{y_{\min }}^{n^{\frac{\beta \gamma }{\left(
\gamma -1\right) }+o_{n}\left( 1\right) }}y^{-\frac{1}{\beta }+\frac{1}{%
\beta \gamma }+o_{y}\left( 1\right) }dy=n^{\frac{\gamma }{\left( \gamma
-1\right) }\left( \beta -1+\frac{1}{\gamma }\right) +o_{n}\left( 1\right) }
\label{ac}
\end{equation}

\textbf{Proof: }We will use the following statement from the appendix in 
\cite{3}: under the assumptions of lemma 2i one has%
\begin{equation}
y^{-1+o_{y}\left( 1\right) }=-\left( D\varphi \right) \circ \varphi
^{-1}\left( y^{-1}\right)  \label{ad}
\end{equation}

Replacing $D\varphi $ by $\varphi $ and $\varphi $ by $F^{\ast }$ we obtain 
\begin{equation}
y^{-1+o_{y}\left( 1\right) }=-\varphi \circ \left( F^{\ast }\right)
^{-1}\left( y^{-1}\right)  \label{aaa}
\end{equation}

Furthermore we will make use of Bernstein inequality in the following form: 
\emph{Let }$X_{1},X_{2},...,X_{n}$\emph{\ be independent random variables
such that }$\left\vert X_{i}\right\vert \leq M$\emph{\ , }$\mathbb{E}\left(
X_{i}\right) =0$\emph{\ and }$b_{n}:=\sum\limits_{i=1}^{n}\mathbb{E}\left(
X_{i}^{2}\right) $\emph{\ . Then for any }$\lambda \geq 0$\emph{\ one has } 
\begin{equation}
\Pr \left( \left\vert \sum\limits_{i=1}^{n}X_{i}\right\vert \geq \lambda
\right) \leq \exp \left( -\frac{\lambda ^{2}}{b_{n}+2M}\right)  \label{22}
\end{equation}

Recall that according to lemma 2 the induced distribution of the \textbf{rv} 
$y=\frac{1}{\varphi ^{\beta }\left( \omega \right) }$ is given by $\psi
\left( y\right) =\frac{1}{y^{\theta +o\left( 1\right) }}$ with $1<\theta $
depending on wether $\varphi \left( \omega \right) $ is of power law type or
not. To prove lemma 3 we therefore have to show 
\begin{equation}
\lim\limits_{n\rightarrow \infty }\frac{\int_{y_{\min }}^{n^{\frac{1}{\theta
-1}}}y\psi \left( y\right) dy}{\frac{1}{n}\sum_{i=1}^{n}y_{i}}=n^{o\left(
1\right) },a.s.
\end{equation}

We will first use Bernstein inequality for the conditioned variable $\tilde{y%
}:=\left( y\mid y<k\right) $ with $k$ large. The distribution of $\tilde{y}$
is given by 
\begin{eqnarray}
\tilde{F}\left( z\right) &=&\Pr \left\{ \tilde{y}<z\right\} =\frac{\Pr
\left\{ y<z\right\} }{\Pr \left\{ y<k\right\} } \\
&=&\frac{\int\limits_{y_{\min }}^{z}y^{-\theta +o\left( 1\right) }dy}{%
\int\limits_{y_{\min }}^{k}y^{-\theta +o\left( 1\right) }dy}=C\left(
k\right) F\left( z\right) ;z\leq k
\end{eqnarray}%
with $C\left( k\right) \rightarrow 1$ as $k\rightarrow \infty $ hence $y$
and $\tilde{y}$ have up to a constant the same distribution. Let $X_{i}:=%
\mathbb{E}\left( \tilde{y}_{i}\right) -\tilde{y}_{i}$ and note that the 
\textbf{rv }$X_{i}$ satisfies the assumptions required in Bernstein
inequality with $M=k$. We have $\mathbb{E}\left( X_{i}^{2}\right) =\mathbb{E}%
\left( \tilde{y} _{i}^{2}\right) -\mathbb{E}^{2}\left( \tilde{y}_{i}\right) $
with $\mathbb{E}\left( \tilde{y}_{i}\right) =\int\limits_{\tilde{y}_{\min
}}^{k}C\left( k\right) \tilde{y}^{1-\theta +o\left( 1\right) }d\tilde{y}%
=k^{2-\theta+o_{k}\left( 1\right) }$ and $\mathbb{E}\left( \tilde{y}%
_{i}^{2}\right) =k^{3-\theta +o_{k}\left( 1\right) }$, hence $\mathbb{E}%
\left( X_{i}^{2}\right) =k^{3-\theta +o_{k}\left( 1\right) }-k^{4-2\theta
+o_{k}\left( 1\right) }$. Since $3-\theta >4-2\theta $ for $\theta >1$ we
have $\mathbb{E}\left( X_{i}^{2}\right) =k^{3-\theta +o_{k}\left( 1\right) }$%
. Furthermore $\sum\limits_{i=1}^{n}\mathbb{E}\left( X_{i}^{2}\right)
+2k=nk^{3-\theta +o_{k}\left( 1\right) }+2k.$ We will take now $k$ as a
function of $n$ namely as the expected maximal value of the \textbf{rv }$y$.
By lemma 1 we then know that the sequence $\left( y_{i}\right) _{1}^{n}$ is
almost surely equal to the sequence $\left( \left( \tilde{y}_{i}\right)
_{1}^{n}\right) $ for $k=n^{\frac{1}{\theta -1}}.$ Let $\lambda
=n^{\varepsilon }$ with $0<\varepsilon $ to be defined later. From equation %
\ref{22} we get 
\begin{eqnarray}
\Pr \left( \left\vert n\mathbb{E}\left( \tilde{y}_{i}\right)
-\sum\limits_{i=1}^{n}y_{i}\right\vert \geq n^{\varepsilon }\right) &\leq
&\exp \left( -\frac{n^{2\varepsilon }}{n\cdot n^{\frac{3-\theta +o_{n}\left(
1\right) }{\theta -1}}+2n^{\frac{1}{\theta -1}}}\right) \\
&\leq &\exp \left( -\frac{n^{2\varepsilon }}{n^{\frac{2+o_{n}\left( 1\right) 
}{\theta -1}}+2n^{\frac{1}{\theta -1}}}\right) \\
&\leq &\exp \left( -n^{2\varepsilon -\frac{2}{\theta -1}+o_{n}\left(
1\right) }\right)
\end{eqnarray}

The right hand side converges to zero for $\varepsilon >\frac{1}{\theta -1}.$
Since furthermore $\mathbb{E}\left( \tilde{y}_{i}\right) =n^{\frac{2-\theta 
}{\theta -1}+o_{n}\left( 1\right) }=n^{\frac{1}{\theta -1}-1+o_{n}\left(
1\right) }$ we obtain $\frac{\mathbb{E}\left( \tilde{y}_{i}\right) }{\frac{1%
}{n}\sum\limits_{i=1}^{n}y_{i}}=\frac{\left( 1+o\left( 1\right) \right)
\int_{y_{\min }}^{n^{\frac{1}{\theta -1}}}y^{1-\theta +o\left( 1\right) }dy}{%
\frac{1}{n}\sum\limits_{i=1}^{n}y_{i}}\leq n^{\varepsilon -\frac{1}{\theta -1%
}+o\left( 1\right) }$ almost surely. Chose a sequence $\left\{ \varepsilon
_{i};\varepsilon _{i}>\frac{1}{\theta -1}\right\} \rightarrow \frac{1}{%
\theta -1}$ For each $\varepsilon _{i}$ we have by Borel-Cantelli lemma $%
\Phi -$ almost surely 
\begin{equation}
\frac{\int_{y_{\min }}^{n^{\frac{1}{\theta -1}}}y^{1-\theta +o\left(
1\right) }dy}{\frac{1}{n}\sum\limits_{i=1}^{n}y_{i}}\leq n^{\varepsilon _{i}-%
\frac{1}{\theta -1}+o_{n}\left( 1\right) }  \label{200}
\end{equation}

Let $G_{i}$ be the set of sequences $\hat{\omega}\in \Sigma $ for which \ref%
{200} holds. Since $\Phi \left( G_{i}\right) =1$ for all $i$ we have for the
countable intersection $G:=\cap _{i}G_{i}$ also $\Phi \left( G\right) =1. $
Clearly for elements in $G$ lemma (3i) is true. The statements in (3 ii) and
(3 iii) follow from a straightforeward computation using \ref{aaa}.

\section{$k$-path connectivity and essential diameter}

Let $d\left( i,j\right) $ be the usual distance between two vertices $i$ and 
$j$ in a given graph $G$. Define as usual the diameter by $diam\left(
G\right) :=\max\limits_{i,j\in V\left( G\right) ,i\neq j}d\left( i,j\right) $
and the component diameter by $diam_{Co}\left( G\right)
:=\max\limits_{Co}diam\left( G\mid _{Co}\right) $ where $G\mid _{Co}$denotes
the restriction of $G$ to a connected component $Co$. We further introduce
the notion of $\varepsilon -$ essential diameter%
\begin{equation}
diam_{ess}^{\left( \varepsilon \right) }\left( G\right)
:=\min\limits_{V^{\ast }\subset V\left( G\right) ;\#V^{\ast }=\left\lceil
\varepsilon \left\vert V\left( G\right) \right\vert \right\rceil }diam\left(
G\mid _{V^{\ast }}\right)
\end{equation}%
The essential diameter is related to the existence of a core ball on the
graph which carries a positive fraction of all the vertices.

The expected path length (\textbf{epl}) $\Delta \left( G\right) $ of a
connected graph $G$ is defined as the mean distance between pairs of
vertices: $\Delta \left( G\right) :=\frac{1}{\binom{n}{2}}%
\sum\limits_{i,j\in V\left( G\right) ;i\neq j}d\left( i,j\right) $. For
non-connected graphs $G$ we define $\Delta \left( G\right) $ as the average
over $d\left( i,j\right) $ where $i,j$ belong to the same connected
component of $G$. It is a well known phenomena that the random variable
describing the diameter or the expected path length of a r.g. space is
typically highly concentrated. Closely related quantities are the
probability $P_{k}$ that there is a path of length $k$ between two at random
chosen vertices and the $k-$ path connectivity $\Gamma _{k}\left(
i,j,G\right) $, the number of paths (without repetition) between $i$ and $j$
in the graph $G$. From $\Gamma _{k}\left( i,j,G\right) $ two natural
quantities can be derived: the maximum $k-$ connectivity $\Gamma _{k}\left(
G\right) :=\max\limits_{i,j}$ $\Gamma _{k}\left( i,j,G\right) $ of a graph $%
G $ and the mean $k-$ connectivity $\bar{\Gamma}_{k}\left( G\right) :=\frac{1%
}{\left\vert E\left( G\right) \right\vert }\sum\limits_{\left( i,j\right)
\in E\left( G\right) }\Gamma _{k}\left( i,j,G\right) $. It turns out that
for large $n$ the r.v. $\Gamma _{k}$ is either close to zero or very large
and the same holds for $\bar{\Gamma}_{k}\left( G\right) $. For random graphs
with a giant component, this "jump" value of $\Gamma _{k}$ is either of the
same order or smaller as $\Delta \left( G\right) $. We are now ready to
state our main theorem:

\textbf{Theorem :}

\textit{i) For} $c>0$\textit{\and } $\alpha \in \left( \frac{1}{2},1\right) $
\textit{there exists} $L\left( c,\alpha \right)$ \textit{such that, with
high probability,\textbf{whp}}, $\Gamma _{k}\left( G\in \mathcal{G}\left( n,%
\hat{\omega}\right) \right)$\textit{is almost 0 for} $k$ \textit{less than} $%
L\left( c,\alpha \right)$ \textit{\and} $\Gamma _{k}\left( G\in \mathcal{G}%
\left( n,\hat{\omega}\right)\right)$ \textit{much bigger than 1 for} $k$ 
\textit{bigger than} $L\left( c,\alpha \right) $ \textit{almost surely in} $%
\hat{\omega}$\textit{.}

\textit{ii) For }$c>\frac{A\left( \alpha \right) }{A\left( \alpha \right) +%
\sqrt{A_{2}}}$\textit{\ and }$\alpha \in \left[ 0,\frac{1}{2}\right) $%
\textit{there exists} $L\left( c,\alpha \right) $ \textit{such that, with
high probability}, $\Gamma _{k}\left(G\in \mathcal{G}\left( n,\hat{\omega}%
\right)\right)$ \textit{is much smaller than 1 for} $k$ \textit{less than} $%
\left( L\left( c,\alpha \right) -1\right) \log n$ \textit{and} $\Gamma
_{k}\left( G\in \mathcal{G}\left( n,\hat{\omega}\right) \right) $ \textit{%
much bigger than 1 for} $k$ \textit{bigger than} $\left( L\left( c,\alpha
\right) +1\right) \log n$ \textit{almost surely in} $\hat{\omega}$\textit{.}

\textit{iii) under the conditions of i) and for }$0<\varepsilon <\varepsilon
_{0}\left( c,\alpha \right) $ \textit{there exists} $L\left( c,\alpha
\right) $ \textit{\ such that, with high probability, }$diam_{ess}^{\left(
\varepsilon \right) }\left( G\right)$ \textit{is less or equal to} $L\left(
c,\alpha \right)$ .

\textbf{Proof:}

i) Let $\alpha \in \left( \frac{1}{2},1\right) $ and $\varphi $ satisfying
the assumptions of lemma 2i. We will start with an estimation of the
expected number $\mathbb{E}_{\mathcal{G}\left( n,\hat{\omega}\right) }\left[
\Gamma _{k}\left( i,j\right) \right] $ of paths of length $k$ between two
vertices $i$ and $j$ for typical $\hat{\omega}$. We will show that there is
a sharp transition value $k_{0}\left( n\right) $ such that, below $k_{0}$
this expectation is close to zero whereas for $k>k_{0}$ it tends to
infinity. Since $\Pr \left\{ d\left( i,j\right) \leq k\mid \mathcal{G}\left(
n,\hat{\omega}\right) \right\} \underset{n\rightarrow \infty }{\rightarrow }%
0 $ for $\mathbb{E}_{\mathcal{G}\left( n,\hat{\omega}\right) }\left[ \Gamma
_{k}\left( i,j\right) \right] \rightarrow 0$ we get as well a lower bound
for the expected path length -\textbf{epl}- and the diameter on $\mathcal{G}%
\left( n,\hat{\omega}\right) $ as $n\rightarrow \infty $.

Let $\Sigma _{ij}^{\left( k\right) }\subset \left\{ i\right\} \times \left(
V_{n}\right) ^{k-1}\times \left\{ j\right\} $ be the set of $k+1-$ strings
of vertices without vertex repetition . Clearly one has $\#\Sigma
_{ij}^{\left( k\right) }=\left( n-2\right) \left( n-3\right) ...\left(
n-k\right) $ for $k>1.$ For $\gamma _{k}=\left(
x_{0}:=i,x_{1},x_{2},...,x_{k-1},x_{k}:=j\right) \in \Sigma _{ij}^{\left(
k\right) }$ we denote by $\mathbf{1}_{\gamma _{k}}\left( G\right) $ the
characteristic function for $\gamma _{k}$ being a path of length $k$ between 
$i$ and $j$, that is all pairs $\left( x_{l},x_{l+1}\right) $ in the string $%
\gamma _{k}$ are edges in $G\in \mathcal{G}\left( n,\hat{\omega}\right) $.
From the model definition we have%
\begin{equation}
\Pr \left\{ \mathbf{1}_{\gamma _{k}}=1\right\}
=\prod\limits_{l=0}^{k-1}p_{x_{l}x_{l+1}}\left( n\right)
\end{equation}%
, which for typical $\hat{\omega}$ and $\alpha \in B\left( \varphi \right) $
takes the form 
\begin{equation}
\Pr \left\{ \mathbf{1}_{\gamma _{k}}=1\right\} =\left[ \frac{c\left(
1+o_{n}\left( 1\right) \right) }{nA\left( \alpha \right) }\right]
^{k}\prod\limits_{l=0}^{k-1}\left( \frac{1}{\varphi ^{\alpha }\left( \omega
_{x_{l}}\right) }+\frac{1}{\varphi ^{\alpha }\left( \omega _{x_{l+1}}\right) 
}\right)  \label{cc}
\end{equation}

With $y_{i}:=\frac{1}{\varphi ^{\alpha }\left( \omega _{i}\right) }$ we get
for the expectation of $\Gamma _{k}\left( i,j\right) $ 
\begin{equation}
\mathbb{E}_{\mathcal{G}\left( n,\hat{\omega}\right) }\left[ \Gamma
_{k}\left( i,j\right) \right] =\sum_{\gamma _{k}\in \Sigma _{ij}^{\left(
k\right) }}\left[ \frac{c\left( 1+o_{n}\left( 1\right) \right) }{nA\left(
\alpha \right) }\right] ^{k}\prod\limits_{l=0}^{k-1}\left(
y_{x_{l}}+y_{x_{l+1}}\right)  \label{c}
\end{equation}%
Multiplying out $\prod\limits_{l=0}^{k-1}\left( y_{x_{l}}+y_{x_{l+1}}\right) 
$ one obtains two different type of terms, namely terms where each $%
y_{x_{l}} $ appears just once and terms where $\left( y_{x_{l}}\right) ^{2}$
appears. It turns out, that the main contribution to \ref{c} comes from
terms with a maximal number of $\left( y_{x_{l}}\right) ^{2}$ involved.

By using the ergodic theorem and lemma 3 we obtain upper and lower bounds on 
$\mathbb{E}_{\mathcal{G}\left( n,\hat{\omega}\right) }\left[ \Gamma
_{k}\left( i,j\right) \right] $ as follows: Since for $k$ fixed we have $%
\Phi $ a.s.%
\begin{equation}
\frac{1}{n\left( n-1\right) ...\left( n-k+1\right) }\sum%
\limits_{x_{1},x_{2},...,x_{k};x_{0}\neq x_{1}\neq x_{2}...\neq
x_{k}}y_{x_{1}}...y_{x_{k}}=\frac{1+o_{n}\left( 1\right) }{n^{k}}%
\sum\limits_{x_{1},x_{2},...,x_{k}}y_{x_{1}}...y_{x_{k}}
\end{equation}%
and%
\begin{equation}
\frac{1+o_{n}\left( 1\right) }{n^{k}}\sum%
\limits_{x_{1},x_{2},...,x_{k}}y_{x_{1}}...y_{x_{k}}=\frac{1+o_{n}\left(
1\right) }{n}\sum\limits_{x_{1}}y_{x_{1}}\cdot \left( \frac{1}{n}%
\sum\limits_{x_{2}}y_{x_{2}}\cdot ....\cdot \left( \frac{1}{n}%
\sum\limits_{x_{k}}y_{x_{k}}\right) \right)
\end{equation}%
applying the ergodic theorem separately to the sums we obtain 
\begin{equation}
\frac{1}{n^{k-1}}\sum\limits_{\gamma _{k}\in \Sigma _{ij}^{\left( k\right)
}}\prod\limits_{x_{l}\in \gamma _{k},0<l<k}y_{x_{l}}=\left( 1+o_{n}\left(
1\right) \right) \left[ \int\limits_{1}^{\infty }\varphi ^{1-\alpha }\left(
\omega \right) d\omega \right] ^{k-1}
\end{equation}

More generally with $\nu =\left( \nu _{1},\nu _{2},...,\nu _{k-1}\right) $
such that $\nu \in \left\{ 0\,;1;2\right\} $ , $\sum \nu _{l}=k-1$ we have 
\begin{equation}
\frac{1}{n^{k-1-g_{0}}}\sum\limits_{\gamma _{k}\in \Sigma _{ij}^{\left(
k\right) }}\prod\limits_{x_{l}\in \gamma _{k}}y_{x_{l}}^{\nu
_{l}}=n^{o_{n}\left( 1\right) }\cdot \left[ \int\limits_{1}^{F^{-1}\left( 1-%
\frac{1}{n}\right) }\varphi ^{1-2\alpha }\left( \omega \right) d\omega %
\right] ^{g_{2}}\cdot \left[ \int\limits_{1}^{\infty }\varphi ^{1-\alpha
}\left( \omega \right) d\omega \right] ^{g_{1}}  \label{d}
\end{equation}%
where $g_{i}:=\#\left\{ \nu _{l}\mid \nu _{l}=i\right\} $ for $i\in \left\{
0\,;1;2\right\} .$ Note that $2\alpha \notin B\left( \varphi \right) $ for $%
\alpha \in \left( \frac{1}{2},1\right) $ and therefore lemma 3 has to be
applied. Using the last expression and the estimation from lemma 3 $%
\int\limits_{1}^{F^{-1}\left( 1-\frac{1}{n}\right) }\varphi ^{1-2\alpha
}\left( \omega \right) d\omega =n^{\delta \left( \alpha \right) }$ where $%
\delta \left( \alpha \right) =2\alpha -1+o_{n}\left( 1\right) ,$ we get the
following bounds:%
\begin{equation}
\mathbb{E}_{\mathcal{G}\left( n,\hat{\omega}\right) }\left[ \Gamma
_{k}\left( i,j\right) \right] \leq \left[ \frac{c\left( 1+o_{n}\left(
1\right) \right) }{A\left( \alpha \right) }\right] ^{k}\frac{n^{\delta
\left( \alpha \right) \left\lfloor k/2\right\rfloor }}{n}A^{k}\left( \alpha
\right) 2^{k}=
\end{equation}%
and 
\begin{equation}
\mathbb{E}_{\mathcal{G}\left( n,\hat{\omega}\right) }\left[ \Gamma
_{k}\left( i,j\right) \right] \geq \left[ \frac{c\left( 1+o_{n}\left(
1\right) \right) }{A\left( \alpha \right) }\right] ^{k}\frac{n^{\delta
\left( \alpha \right) \left\lfloor k/2\right\rfloor }}{n}
\end{equation}

Clearly $\mathbb{E}_{\mathcal{G}\left( n,\hat{\omega}\right) }\left[ \Gamma
_{k}\left( i,j\right) \right] \rightarrow 0$ for $k<\frac{2}{\delta \left(
\alpha \right) }$ and $\mathbb{E}_{\mathcal{G}\left( n,\hat{\omega}\right) }%
\left[ \Gamma _{k}\left( i,j\right) \right] \rightarrow \infty $ for $k>%
\frac{2}{\delta \left( \alpha \right) }+1.$ To finish the proof of part i)
of theorem 1 we still have to establish a relation between $\mathbb{E}_{%
\mathcal{G}\left( n,\hat{\omega}\right) }\left[ \Gamma _{k}\left( i,j\right) %
\right] $ and $\mathbb{\Gamma }_{k}\left( G\right), $ respectively $\bar{%
\Gamma}_{k}\left( G\right) $. For the case when $\varphi $ is itself a power
law distribution $\frac{const}{\omega ^{\gamma }}$ one gets in complete
analogy to the above computations the following estimation.

ii) Let $\alpha \in \left( 0,\frac{1}{2}\right) $. The first part of the
argumentation in the previous section (till formula \ref{c}) remains
unchanged since it was independent of the $\alpha -$ value. In formula \ref%
{d} we have for $\alpha <\frac{1}{2}$ and $\varphi $ decaying faster than
any power law the estimation%
\begin{equation}
\int\limits_{1}^{F^{-1}\left( 1-\frac{1}{n}\right) }\varphi ^{1-2\alpha
}\left( \omega \right) d\omega <\int\limits_{1}^{\infty }\varphi ^{1-2\alpha
}\left( \omega \right) d\omega =:A_{2}>A:=A\left( \alpha \right)
\end{equation}

This gives in connection with formula \ref{c}:%
\begin{equation}
\mathbb{E}_{\mathcal{G}\left( n,\hat{\omega}\right) }\left[ \Gamma
_{k}\left( i,j\right) \right] \leq \left[ \frac{c\left( 1+o_{n}\left(
1\right) \right) }{A\left( \alpha \right) }\right] ^{k}\frac{A_{3}}{n}
\end{equation}%
and 
\begin{equation}
\mathbb{E}_{\mathcal{G}\left( n,\hat{\omega}\right) }\left[ \Gamma
_{k}\left( i,j\right) \right] \geq \left[ \frac{c\left( 1+o_{n}\left(
1\right) \right) }{A\left( \alpha \right) }\right] ^{k}\frac{A_{4}}{n}
\end{equation}%
where the constants $A_{3}$ and $A_{4}$ are given by $A_{3}=\left(
2A_{2}\right) ^{k/2}$ and $A_{4}=A^{k}$. Hence for $c\ $sufficiently small
and $\Phi $ a.s. \textbf{whp }$G\in \mathcal{G}\left( n,\hat{\omega}\right) $
has no giant component since otherwise the expected $k-$ path number has to
increase with $k$. To get better bounds on the the jump value for the $k-$
path number for $c>1$ we have to estimate the sum over the products $%
\prod\limits_{l=0}^{k-1}\left( y_{x_{l}}+y_{x_{l+1}}\right) $ a bit more
precisely. Let $C_{m}^{\left( k\right) }$ be the number of terms $a_{j}$ in
the product $\prod\limits_{l=0}^{k-1}\left( y_{x_{l}}+y_{x_{l+1}}\right)
=\sum a_{j}$ with $a_{j}=\prod y_{x_{l}}^{\nu _{l}}$ for some vector $\nu
=\left( \nu _{1},\nu _{2},...,\nu _{k-1}\right) $, $\nu _{i}\in \left\{
0\,;1;2\right\} ,$ where exactly for $m$ values of $i$ the exponent $\nu
_{i}=2$ appears. Replacing the sums over the different $x_{l}$ for $\Phi -$
typical $\hat{\omega}$ by the corresponding integral we get 
\begin{equation}
\mathbb{E}_{\mathcal{G}\left( n,\hat{\omega}\right) }\left[ \Gamma
_{k}\left( i,j\right) \right] =\frac{1}{n}\left[ \frac{c\left( 1+o_{n}\left(
1\right) \right) }{A\left( \alpha \right) }\right] ^{k}\sum\limits_{m=0}^{m=%
\left\lfloor \frac{k}{2}\right\rfloor }C_{m}^{\left( k\right)
}A_{2}^{m}A^{k-2m}  \label{100}
\end{equation}

The coefficients $C_{m}^{\left( k\right) }$ can nicely be interpreted in
terms of a topological Markov chain. Consider the $4-$ state Markov chain
with transition matrix 
\begin{equation}
\left( a_{ij}\right) =%
\begin{pmatrix}
1 & 0 & 1 & 0 \\ 
0 & 1 & 0 & 1 \\ 
0 & 1 & 0 & 1 \\ 
1 & 0 & 1 & 0%
\end{pmatrix}%
\end{equation}

The states correspond to the different types of product terms in $\left(
y_{x_{l}}+y_{x_{l+1}}\right) \left( y_{x_{l+1}}+y_{x_{l+2}}\right) $, namely 
$1\doteq y_{x_{l}}y_{x_{l+1}},$ $2\doteq y_{x_{l+1}}y_{x_{l+2}},3\doteq
y_{x_{l}}y_{x_{l+2}}$ and $4\doteq y_{x_{l+1}}y_{x_{l+1}}$. Then $%
C_{m}^{\left( k\right) }$ is the number of words of length $k$ with exactly $%
m-$ times the symbol $4$ appearing.

We illustrate the construction by an example. Let $k=4$ and consider the
products in the path probability sum 
\begin{eqnarray}
\prod\limits_{l=0}^{3}\left( y_{x_{l}}+y_{x_{l+1}}\right)
&=&(y_{x_{0}}y_{x_{1}}y_{x_{2}}y_{x_{3}}+y_{x_{0}}y_{x_{1}}y_{x_{2}}y_{x_{4}}+y_{x_{0}}y_{x_{1}}\allowbreak y_{x_{3}}y_{x_{4}}+
\notag \\
&&+y_{x_{0}}y_{x_{2}}y_{x_{3}}y_{x_{4}}+y_{x_{1}}y_{x_{2}}y_{x_{3}}y_{x_{4}}+\allowbreak y_{x_{0}}y_{x_{1}}y_{x_{3}}^{2}+
\notag \\
&&+y_{x_{0}}y_{x_{2}}y_{x_{3}}^{2}+y_{x_{0}}y_{x_{2}}^{2}y_{x_{3}}+%
\allowbreak y_{x_{0}}y_{x_{2}}^{2}y_{x_{4}}+ \\
&&+y_{x_{1}}y_{x_{2}}y_{x_{3}}^{2}+y_{x_{1}}y_{x_{2}}^{2}y_{x_{3}}+%
\allowbreak y_{x_{1}}^{2}y_{x_{2}}y_{x_{3}}+y_{x_{1}}y_{x_{2}}^{2}y_{x_{4}} 
\notag \\
&&+y_{x_{1}}^{2}y_{x_{2}}y_{x_{4}}+\allowbreak
y_{x_{1}}^{2}y_{x_{3}}y_{x_{4}}+y_{x_{1}}^{2}y_{x_{3}}^{2})  \notag
\end{eqnarray}

Each product corresponds uniquely to one symbolic word e.g. 
\begin{eqnarray}
y_{x_{0}}y_{x_{1}}y_{x_{2}}y_{x_{3}} &\doteq &111 \\
y_{x_{0}}y_{x_{1}}\allowbreak y_{x_{3}}y_{x_{4}} &\doteq &132 \\
y_{x_{1}}y_{x_{2}}y_{x_{3}}y_{x_{4}} &\doteq &222 \\
y_{x_{1}}y_{x_{2}}^{2}y_{x_{3}} &\doteq &241 \\
y_{x_{1}}^{2}y_{x_{3}}^{2} &\doteq &434
\end{eqnarray}%
and so on.

To estimate $S\left( k,A_{2},A\right) :=$ $\sum\limits_{m=0}^{m=\left\lfloor 
\frac{k}{2}\right\rfloor }C_{m}^{\left( k\right) }A_{2}^{m}A^{k-2m}$ we will
first derive a recursion relation for the numbers $C_{m}^{\left( k\right) }$
conditioned on the end-value of the $k-$string. Let therefore $%
C_{m,a}^{\left( k\right) }$ be the number of strings of length $k$ ending in
state $a$ , ( $a=1,2,3,4$ ) with exactly $m$ times the symbol $4$ appearing.
We get the following recursion: 
\begin{eqnarray}
C_{m,1}^{\left( k+1\right) } &=&C_{m,1}^{\left( k\right) }+C_{m,4}^{\left(
k\right) } \\
C_{m,2}^{\left( k+1\right) } &=&C_{m,2}^{\left( k\right) }+C_{m,3}^{\left(
k\right) } \\
C_{m,3}^{\left( k+1\right) } &=&C_{m,1}^{\left( k\right) }+C_{m,4}^{\left(
k\right) } \\
C_{m,4}^{\left( k+1\right) } &=&C_{m-1,2}^{\left( k\right)
}+C_{m-1,3}^{\left( k\right) }
\end{eqnarray}%
for $k\geq 1$ and $m>0$. The initial conditions are: $C_{0,a}^{\left(
1\right) }=1$ for $a=1,2,3$, $C_{0,4}^{\left( 1\right) }=0,$ $%
C_{1,a}^{\left( 1\right) }=0$ for $a=1,2,3,C_{1,4}^{1}=1.$ With $%
X_{m}^{k}:=C_{m,1}^{\left( k\right) }+C_{m,4}^{\left( k\right) }$ and $%
Y_{m}^{k}:=C_{m,2}^{\left( k\right) }+C_{m,3}^{\left( k\right) }$ this
reduces to 
\begin{eqnarray}
X_{m}^{\left( k+1\right) } &=&X_{m}^{\left( k\right) }+Y_{m-1}^{\left(
k\right) } \\
Y_{m}^{\left( k+1\right) } &=&X_{m}^{\left( k\right) }+Y_{m}^{\left(
k\right) }
\end{eqnarray}

Using the generating functions $f_{X}\left( k,z\right)
:=\sum\limits_{m}X_{m}^{\left( k\right) }z^{m}$ and $f_{Y}\left( k,z\right)
:=\sum\limits_{m}Y_{m}^{\left( k\right) }z^{m}$ the above recursion relation
for the $C_{m,a}^{\left( k\right) }$ can be translated into a recursion
relation for the generating functions:%
\begin{equation}
\begin{pmatrix}
f_{X}\left( k+1,z\right) \\ 
f_{Y}\left( k+1,z\right)%
\end{pmatrix}%
=%
\begin{pmatrix}
1 & z \\ 
1 & 1%
\end{pmatrix}%
\begin{pmatrix}
f_{X}\left( k,z\right) \\ 
f_{Y}\left( k,z\right)%
\end{pmatrix}%
,k\geq 1
\end{equation}%
with initial conditions $f_{X}\left( 1,z\right) =1+z$ and $f_{Y}\left(
1,z\right) =2.$ The transition matrix has eigenvalues $\left\{ \lambda _{1}=%
\sqrt{z}+1,\lambda _{2}=1-\sqrt{z}\right\} $ and therefore we obtain%
\begin{eqnarray}
f_{X}\left( k,z\right) &=&\left( \frac{1}{2}\left( 1+z\right) +\sqrt{z}%
\right) \lambda _{1}^{k-1}+\left( \frac{1}{2}\left( 1+z\right) -\sqrt{z}%
\right) \lambda _{2}^{k-1} \\
f_{Y}\left( k,z\right) &=&\left( 1+\frac{1}{2\sqrt{z}}+\frac{\sqrt{z}}{2}%
\right) \lambda _{1}^{k-1}+\left( 1-\frac{1}{2\sqrt{z}}-\frac{\sqrt{z}}{2}%
\right) \lambda _{2}^{k-1}
\end{eqnarray}%
for $k\geq 1$. For $z=\frac{A_{2}}{A^{2}}$ and since $S\left(
k,A_{2},A\right) :=$ $A^{k}\sum\limits_{m=0}^{m=\left\lfloor \frac{k}{2}%
\right\rfloor }C_{m}^{\left( k\right) }\left( \frac{A_{2}}{A^{2}}\right)
^{m} $ we get $S\left( k,A_{2},A\right) =const\cdot \left( A+\sqrt{A_{2}}%
+o_{k}\left( 1\right) \right) ^{k}.$ Inserting this into equation \ref{100}
we obtain%
\begin{eqnarray}
\mathbb{E}_{\mathcal{G}\left( n,\hat{\omega}\right) }\left[ \Gamma
_{k}\left( i,j\right) \right] &=&\frac{const}{n}\left[ \frac{c\left(
1+o_{n}\left( 1\right) \right) }{A\left( \alpha \right) }\right] ^{k}\left(
A\left( \alpha \right) +\sqrt{A_{2}}+o_{k}\left( 1\right) \right) ^{k} 
\notag \\
&=&\frac{const}{n}\left( 1+o_{n}\left( 1\right) \right) ^{k}\left( c+\frac{c%
\sqrt{A_{2}}}{A\left( \alpha \right) }\right) ^{k}
\end{eqnarray}

For $c+\frac{c\sqrt{A_{2}}}{A\left( \alpha \right) }<1$ the expected $k-$%
path number converges to zero for each $k$ and hence the corresponding
random graph space has \textbf{whp} no giant component. For $B:=c+\frac{c%
\sqrt{A_{2}}}{A\left( \alpha \right) }>1$ and $k=\frac{\left( 1+\varepsilon
\right) \log n}{\log B}$ , $\varepsilon >0$ the expected $k$ $-$ path number
goes to infinity.

It remains to prove that the above phase transition happens \textbf{whp} in $%
G\in \mathcal{G}\left( n,\hat{\omega}\right) $ for typical $\hat{\omega}.$ A
standard technique to obtain such results is the so called first and second
moment method. From Tchebychev inequality one has for a discrete positive
random variable $X$ with $\mathbb{E}\left( X\right) =a$ 
\begin{equation}
\Pr \left\{ X=0\right\} \leq \frac{\mathbb{E}\left( X^{2}\right) }{a^{2}}-1
\end{equation}%
. By the Markov inequality we have further 
\begin{equation}
\Pr \left\{ X\geq t\right\} \leq \frac{a}{t}
\end{equation}%
. We will first show that 
\begin{equation}
\frac{\mathbb{E}_{\mathcal{G}\left( n,\hat{\omega}\right) }\left( \left[
\Gamma _{k}\left( i,j\right) \right] ^{2}\right) }{\left[ \mathbb{E}_{%
\mathcal{G}\left( n,\hat{\omega}\right) }\left( \left[ \Gamma _{k}\left(
i,j\right) \right] \right) \right] ^{2}}\underset{n\rightarrow \infty }{%
\rightarrow }1  \label{1000}
\end{equation}%
for almost every $\hat{\omega}$ and $k>k_{c}$ where $k_{c}$ is the phase
transition value in Theorem 1. For a given string $\gamma _{k}=\left\{
i=x_{0},x_{1},...,x_{k}=j\right\} \in \Sigma _{ij}^{\left( k\right) }$ let%
\begin{equation}
B_{l}\left( \gamma _{k}\right) :=\left\{ \gamma _{k}^{\prime }\in \Sigma
_{ij}^{\left( k\right) }:\gamma _{k}^{\prime }\text{ has }l\text{ pairs }%
\left( x_{m,}x_{m+1}\right) \text{ in common with }\gamma _{k}\right\}
\end{equation}%
and 
\begin{equation}
A_{l}\left( \gamma _{k}\right) :=\left\{ \gamma _{k}^{\prime }\in \Sigma
_{ij}^{\left( k\right) }:\gamma _{k}^{\prime }\text{ has }l\text{ interior
vertices }x_{m}\text{ in common with }\gamma _{k}\right\}
\end{equation}%
where the the interior vertices are all $x_{m}$ with $1\leq m\leq k-1.$ The
index $in\left( \gamma _{k}^{\prime },\gamma _{k}\right) $ of a string $%
\gamma _{k}^{\prime }$ with respect to a given $\gamma _{k}$ is defined as
the unique integer pair $\left( I^{A},I^{B}\right) $ such that $\gamma
_{k}^{\prime }\in A_{I^{A}}\left( \gamma _{k}\right) \cap B_{I^{B}}\left(
\gamma _{k}\right) $. Since for any fixed $G$ we have by definition $\Gamma
_{k}\left( i,j\right) =\sum\limits_{\gamma \in \Sigma _{ij}^{\left( k\right)
}}\mathbf{1}_{\gamma }\left( G\right) $ we get%
\begin{eqnarray}
\left[ \Gamma _{k}\left( i,j\right) \right] ^{2} &=&\sum\limits_{\gamma \in
\Sigma _{ij}^{\left( k\right) }}\mathbf{1}_{\gamma }\left( G\right)
\sum\limits_{\gamma ^{\prime }\in \Sigma _{ij}^{\left( k\right) }}\mathbf{1}%
_{\gamma ^{\prime }}\left( G\right) \\
&=&\sum\limits_{\gamma \in \Sigma _{ij}^{\left( k\right) }}\mathbf{1}%
_{\gamma }\left( G\right) \sum\limits_{l,m}\sum\limits_{\gamma ^{\prime
}:in\left( \gamma _{k}^{\prime },\gamma _{k}\right) =\left( l,m\right) }%
\mathbf{1}_{\gamma ^{\prime }}\left( G\right)
\end{eqnarray}%
. We will show that the main contribution to the above sum is due to strings
with index $\left( 0,0\right) ,$ that is for independent pairs of strings.
More precisely the following proposition holds:

\begin{equation}
\mathbb{E}_{\mathcal{G}\left( n,\hat{\omega}\right) }\left( \left[ \Gamma
_{k}\left( i,j\right) \right] ^{2}\right) =\mathbb{E}_{\mathcal{G}\left( n,%
\hat{\omega}\right) }\left( \sum\limits_{\gamma \in \Sigma _{ij}^{\left(
k\right) }}\mathbf{1}_{\gamma }\left( G\right) \sum\limits_{\gamma ^{\prime
}:in\left( \gamma _{k}^{\prime },\gamma _{k}\right) =\left( 0,0\right) }%
\mathbf{1}_{\gamma ^{\prime }}\left( G\right) \right) +O\left( n^{-1}\right) 
\label{xxx}
\end{equation}%
which implies that for $\mathbb{E}_{\mathcal{G}\left( n,\hat{\omega}\right)
}\left( \sum\limits_{\gamma \in \Sigma _{ij}^{\left( k\right) }}\mathbf{1}%
_{\gamma }\left( G\right) \sum\limits_{\gamma ^{\prime }:in\left( \gamma
_{k}^{\prime },\gamma _{k}\right) =\left( 0,0\right) }\mathbf{1}_{\gamma
^{\prime }}\left( G\right) \right) >1$ . To see this we need an estimation
of the cardinality of $A_{l}\left( \gamma _{k}\right) $ and of the
expectation $\mathbb{E}_{\mathcal{G}\left( n,\hat{\omega}\right) }\left( 
\mathbf{1}_{\gamma _{k}}\cdot \mathbf{1}_{\gamma _{k}^{\prime }}\right) $
for given index of $\gamma _{k}^{\prime }.$ Clearly for $in\left( \gamma
_{k}^{\prime },\gamma _{k}\right) =\left( 0,0\right) $ one has $\mathbb{E}_{%
\mathcal{G}\left( n,\hat{\omega}\right) }\left( \mathbf{1}_{\gamma
_{k}}\cdot \mathbf{1}_{\gamma _{k}^{\prime }}\right) =\mathbb{E}_{\mathcal{G}%
\left( n,\hat{\omega}\right) }\left( \mathbf{1}_{\gamma _{k}}\right) \cdot 
\mathbb{E}_{\mathcal{G}\left( n,\hat{\omega}\right) }\left( \mathbf{1}%
_{\gamma _{k}^{\prime }}\right) $ since the strings are independend.
Furthermore we have $\#A_{0}\left( \gamma _{k}\right)
=\prod\limits_{l=1}^{k-1}\left( n-k-l\right) =n^{k-1}\left( 1+o_{n}\left(
1\right) \right) $ and $\#\Sigma _{ij}^{\left( k\right) }=n^{k-1}\left(
1+o_{n}\left( 1\right) \right) .$ Since the cardinalities of $A_{0}\left(
\gamma _{k}\right) $ and $\Sigma _{ij}^{\left( k\right) }$are of the same
order we obtain 
\begin{eqnarray}
&&\mathbb{E}_{\mathcal{G}\left( n,\hat{\omega}\right) }\left(
\sum\limits_{\gamma \in \Sigma _{ij}^{\left( k\right) }}\mathbf{1}_{\gamma
}\left( G\right) \sum\limits_{\gamma ^{\prime }:in\left( \gamma _{k}^{\prime
},\gamma _{k}\right) =\left( 0,0\right) }\mathbf{1}_{\gamma ^{\prime
}}\left( G\right) \right)  \\
&=&\sum\limits_{\gamma \in \Sigma _{ij}^{\left( k\right) }}\mathbb{E}_{%
\mathcal{G}\left( n,\hat{\omega}\right) }\left( \mathbf{1}_{\gamma }\left(
G\right) \right) \sum\limits_{\gamma ^{\prime }:in\left( \gamma _{k}^{\prime
},\gamma _{k}\right) =\left( 0,0\right) }\mathbb{E}_{\mathcal{G}\left( n,%
\hat{\omega}\right) }\left( \mathbf{1}_{\gamma ^{\prime }}\left( G\right)
\right)   \notag \\
&=&\left( 1+o_{n}\left( 1\right) \right) \sum\limits_{\gamma \in \Sigma
_{ij}^{\left( k\right) }}\mathbb{E}_{\mathcal{G}\left( n,\hat{\omega}\right)
}\left( \mathbf{1}_{\gamma }\left( G\right) \right) \sum\limits_{\gamma
^{\prime }\in \Sigma _{ij}^{\left( k\right) }}\mathbb{E}_{\mathcal{G}\left(
n,\hat{\omega}\right) }\left( \mathbf{1}_{\gamma ^{\prime }}\left( G\right)
\right)   \notag \\
&=&\left( 1+o_{n}\left( 1\right) \right) \left[ \mathbb{E}_{\mathcal{G}%
\left( n,\hat{\omega}\right) }\left( \Gamma _{k}\left( i,j\right) \right) %
\right] ^{2}  \notag
\end{eqnarray}%
. We will now estimate the contribution of the summation over indices with $%
in\left( \gamma _{k}^{\prime },\gamma _{k}\right) =\left( a,b\right) \neq
\left( 0,0\right) .$ First note that for a pair of strings $\left(
i_{0},j_{0}\right) ,\left( i_{1},j_{0}\right) $ where $i_{0}\neq i_{1}$, we
have 
\begin{equation}
\mathbb{E}_{\mathcal{G}\left( n,\hat{\omega}\right) }\left( \mathbf{1}%
_{\left( i_{0},j_{0}\right) }\cdot \mathbf{1}_{\left( i_{1},j_{0}\right)
}\right) =\left[ \frac{c\left( 1+o_{n}\left( 1\right) \right) }{nA\left(
\alpha \right) }\right] ^{2}\left( y_{i_{0}}+y_{j_{0}}\right) \left(
y_{i_{1}}+y_{j_{0}}\right) 
\end{equation}%
. Furthermore we have trivially for identical pairs $\left(
i_{0},j_{0}\right) ,\left( i_{0},j_{0}\right) $%
\begin{equation}
\mathbb{E}_{\mathcal{G}\left( n,\hat{\omega}\right) }\left( \mathbf{1}%
_{\left( i_{0},j_{0}\right) }\cdot \mathbf{1}_{\left( i_{0},j_{0}\right)
}\right) =\mathbb{E}_{\mathcal{G}\left( n,\hat{\omega}\right) }\left( 
\mathbf{1}_{\left( i_{0},j_{0}\right) }\right) 
\end{equation}%
and for pairs $\left( i_{0},j_{0}\right) ,\left( i_{1},j_{1}\right) $ where $%
i_{0}\neq i_{1}$ and $j_{0}\neq j_{1}$ (independent pairs)%
\begin{equation}
\mathbb{E}_{\mathcal{G}\left( n,\hat{\omega}\right) }\left( \mathbf{1}%
_{\left( i_{0},j_{0}\right) }\cdot \mathbf{1}_{\left( i_{1},j_{1}\right)
}\right) =\mathbb{E}_{\mathcal{G}\left( n,\hat{\omega}\right) }\left( 
\mathbf{1}_{\left( i_{0},j_{0}\right) }\right) \cdot \mathbb{E}_{\mathcal{G}%
\left( n,\hat{\omega}\right) }\left( \mathbf{1}_{\left( i_{1},j_{1}\right)
}\right) 
\end{equation}%
. The cardinality of $A_{l}\left( \gamma _{k}\right) $ can easily be
estimated by $\#A_{l}\left( \gamma _{k}\right) =\binom{n-k-1}{k-1-l}\simeq
n^{k-1-l}\left( 1+o_{n}\left( 1\right) \right) .$ This provides also upper
bounds on the cardinality of $B_{l}\left( \gamma _{k}\right) $ since $%
I^{B}\leq I^{A}$ always holds. To see that for $\left( l,m\right) >\left(
0,0\right) $ one has for $\gamma \in \Sigma _{ij}^{\left( k\right) }$

\begin{equation}
\mathbb{E}_{\mathcal{G}\left( n,\hat{\omega}\right) }\left( \mathbf{1}%
_{\gamma }\left( G\right) \sum\limits_{\gamma ^{\prime }:in\left( \gamma
_{k}^{\prime },\gamma \right) =\left( l,m\right) }\mathbf{1}_{\gamma
^{\prime }}\left( G\right) \right) =o\left( \mathbb{E}_{\mathcal{G}\left( n,%
\hat{\omega}\right) }\left( \mathbf{1}_{\gamma }\left( G\right)
\sum\limits_{\gamma ^{\prime }:in\left( \gamma _{k}^{\prime },\gamma \right)
=\left( 0,0\right) }\mathbf{1}_{\gamma ^{\prime }}\left( G\right) \right)
\right)   \label{111}
\end{equation}%
observe first that  the nominator  in the above expression is given by $n$
to the power of the number of different variables about which the ergodic
averages take place- that is $k-1+I^{A}.$ Second the nominator is given by $n
$ to the power of the number of different edges involved in the paths $%
\gamma $ and $\gamma ^{\prime }$ , that is $k+m$ and edges on the left hand
side and $2k$ edges in the right hand side. Therefore for positive index we
can conclude \cite{111} from which the claims of the theorem follow
straightforward.

\section{Comments and conclusions}

In this article we have discussed how the $k-$ path connectivity and the
essential diameter depend on the parameter $\alpha $ expressing for
Cameo-graphs the affinity of the vertices for the property encoded in the
random variable $\omega $. We obtained rigorous results for the
corresponding phase transitions. In a forthcoming paper we will analyze the
component size and the diameter of the largest component in Cameo-graphs.
For this the theory developed in \cite{9} seems to be appropriate.

Cameo-graphs define a model of inhomogeneous random graphs. There remain
many interesting and challenging mathematical questions for further
research. For example both $\omega $ and $\alpha $ are properties of the
vertices, so in applications they take in general different values for
different vertices. We consider in this paper only $\omega $ as a random
variable. What we need also is a knowledge on the distribution of the
affinity $\alpha $. See e.g. \cite{7} for first results in this direction.

The study of spreading processes on complex random networks is a topic of
great interest in many fields and has received considerable attention during
the last years. Epidemiology modelling can be used in planning and
evaluating various prediction scenarios. See \cite{8} for an application to
the study of corruption as an epidemic process.

A natural alternative to the choice of the pairing probabilities given in
equation \ref{1} is a multiplicative kernel function. Comparable results for
this case are easy to obtain by using similar arguments as in this article.
Some parts of the combinatorics become even considerably simpler.

\bigskip

\textbf{Acknowledgements: }We would like to thank for the support of the
Volkswagen Foundation and of the DFG-research group 399 "Spectral Analysis,
Asymptotic Distributions and Stochastic Dynamics" and also of the Center for
Theoretical Physics in Marseille.

\section{Appendix}

We prove here the technical statement used in the proof of lemma 3, namely
that under the assumptions of lemma 2i the following holds:

\begin{equation}
\frac{-\varphi \circ \varphi ^{-1}\left( y^{\frac{-1}{\alpha }}\right) }{%
\left( D\varphi \right) \left( \varphi ^{-1}\left( y^{\frac{-1}{\alpha }%
}\right) \right) }=y^{o_{y}\left( 1\right) }
\end{equation}%
Since $\varphi \left( \omega \right) $ decays faster then any power law, we
have 
\begin{equation}
\varphi \left( \omega \right) <\frac{1}{\omega ^{l}}\text{ for any }l\text{
and }\omega >\omega _{0}\left( l\right) \ .
\end{equation}%
Since $\varphi ^{-1}\left( \frac{1}{y^{\frac{1}{\alpha }}}\right) $ goes to
infinity for $y\rightarrow \infty $ we have to show 
\begin{equation}
\frac{-\varphi \left( x\right) }{D\varphi \left( x\right) }=\left[ \varphi
\left( x\right) \right] ^{o_{x}\left( 1\right) }\ .  \label{A}
\end{equation}%
The last formula states that the negative logarithmic derivative of $\varphi 
$ should not become to large or to small compared to $\varphi $ respectively 
$\frac{1}{\varphi }$ . For the following it is convenient to set $\varphi
\left( x\right) =e^{-g\left( x\right) }$ with $g\left( x\right) \rightarrow
\infty $ and rewrite formula \ref{A} as 
\begin{equation}
e^{-\mu g\left( x\right) }<\frac{1}{Dg\left( x\right) }<e^{\mu g\left(
x\right) }\text{ for }\mu \in \left( 0,\mu _{0}>0\right) \text{ and }%
x>x_{0}\left( \mu \right) \ .  \label{B}
\end{equation}%
Assume that formula \ref{B} is not true with respect to the right hand side.
Then we have for a sequence of values $\left\{ x_{i}\right\} $ and open
intervals $I_{i}$ around the $x_{i}$ and some function $a\left( x\right) $ 
\begin{equation}
\frac{1}{Dg\left( x\right) }=e^{\mu g\left( x\right) }a\left( x\right) \text{
and }a\left( x\right) >1\text{ for }x\in I_{i}\ .
\end{equation}%
Integrating the last equation gives 
\begin{equation}
e^{\mu g\left( x\right) }=e^{\mu g\left( x_{0}\right) }+\mu
\int\limits_{x_{0}}^{x}\frac{1}{a\left( z\right) }dz\ .
\end{equation}%
Since our assumption on $D\left[ \frac{1}{\varphi \left( \omega \right) }%
\right] ^{\mu }$ to be monotonous for $\mu >0$ and $\omega >\omega
_{0}\left( \mu \right) $ implies $a\left( x\right) >1$ eventually we
conclude that 
\begin{equation}
e^{\mu g\left( x\right) }<e^{\mu g\left( x_{0}\right) }+\mu \left(
x-x_{0}\right) \ .  \label{C}
\end{equation}%
But the fast decay condition for $\varphi \left( x\right) $ expresses a
growth condition for $g\left( x\right) $ namely for all $k$ 
\begin{equation}
g\left( x\right) >k\log x;\text{ }x>x_{0}\left( k\right)
\end{equation}%
which clearly contradicts formula \ref{C}. It remains to show that the left
hand side of formula \ref{B} also holds. Assuming the converse we get 
\begin{equation}
\frac{1}{Dg\left( x\right) }=e^{-\mu g\left( x\right) }\frac{1}{a\left(
x\right) }\text{ and }a\left( x\right) >1\text{ for }x\in I_{i}
\end{equation}%
and after integration 
\begin{equation}
e^{-\mu g\left( x\right) }=e^{-\mu g\left( x_{0}\right) }-\mu
\int\limits_{x_{0}}^{x}a\left( z\right) dz\ .
\end{equation}%
The monotonicity condition again implies $a\left( x\right) >1$ eventually,
hence 
\begin{equation}
e^{-\mu g\left( x\right) }<e^{-\mu g\left( x_{0}\right) }-\mu \left(
x-x_{0}\right)
\end{equation}%
and a clear contradiction since the right hand side becomes negative for
large values of $x$.$\square $

\end{document}